\documentclass{groups13}
\usepackage{curves,url}
\newtheorem{conj}[thm]{Conjecture}

\newcommand{\rank}{\mathop{\mathrm{rank}}}
\newcommand{\sym}{S_n}
\newcommand{\trans}{T_n}
\newcommand{\Sym}{\mathop{\mathrm{Sym}}}
\newcommand{\End}{\mathop{\mathrm{End}}}
\newcommand{\Aut}{\mathop{\mathrm{Aut}}}
\newcommand{\Gr}{\mathop{\mathrm{Gr}}}
\newcommand{\pgl}{\mathop{\mathrm{PGL}}}
\newcommand{\agl}{\mathop{\mathrm{AGL}}}
\newcommand{\gl}{\mathop{\mathrm{GL}}}
\newcommand{\GF}{\mathop{\mathrm{GF}}}
\newcommand{\psl}{\mathop{\mathrm{PSL}}}
\newcommand{\Sz}{\mathop{\mathrm{Sz}}}
\newcommand{\pgaml}{\mathop{\mathrm{P\Gamma L}}}

\begin{document}

\maketitle%
{PERMUTATION GROUPS AND TRANSFORMATION SEMIGROUPS: RESULTS AND PROBLEMS}%
{JO\~AO ARA\'UJO$^{\ast}$ and PETER J. CAMERON$^{\dagger}$}%
{$^{\ast}$Universidade Aberta and Centro de Algebra, Universidade de Lisboa, Av.\
Gama Pinto 2, 1649-003 Lisboa, Portugal\newline
Email: jaraujo@ptmat.fc.ul.pt\\[2pt]
$^{\dagger}$Mathematical Institute, University of St Andrews, North
Haugh, St Andrews, Fife, KY16 9SS, U.K.\newline
Email: pjc@mcs.st-andrews.ac.uk}

\begin{abstract}
J.M. Howie, the influential St Andrews semigroupist, claimed that we
value an area of pure mathematics to the extent that (a) it gives rise to
arguments that are deep and elegant, and (b) it has interesting
interconnections with other parts of pure mathematics.

This paper surveys some recent results on the transformation semigroup
generated by a permutation group $G$ and a single non-permutation $a$. Our
particular concern is the influence that properties of $G$ (related to
homogeneity, transitivity and primitivity) have on the structure of the
semigroup. In the first part of the paper, we consider properties of
$S=\langle G,a\rangle$ such as regularity and idempotent generation. The second
is a brief report on the synchronization project, which aims to decide
in what circumstances $S$ contains an element of rank~$1$. The paper closes
with a list of open problems on permutation groups and linear groups, and some
comments about the impact on semigroups are provided.

These two research directions outlined above lead to very interesting and
challenging problems on primitive permutation groups whose solutions require
combining results from several different areas of mathematics, certainly
fulfilling both of Howie's elegance and value tests in a new and fascinating
way.
\end{abstract}

\section{Regularity and generation}

\subsection{Introduction}

How can group theory help the study of semigroups?

If a semigroup has a large group of units, we can apply group theory to it.
But there may not be any units at all! According to a widespread  belief,
almost all finite semigroups have only one idempotent, which is a zero, not an identity 
(see \cite{KRS} and \cite{JMS}). This conjecture, however, should not deter
us from the general goal of investigating how the group of units shapes the
structure of the semigroup. Infinitely many families of finite semigroups, and the most interesting,  are composed by semigroups
with group of units. Some of those families are
interesting enough to keep many mathematicians busy their entire lives; in
fact a unique family of finite semigroups, the endomorphism semigroups of
vector spaces over finite fields, has been keeping experts in linear algebra
busy for more than a century.

Regarding the general question of how the group of units can shape the
structure of the semigroup, an especially promising area is the theory of
\emph{transformation semigroups}, that is, semigroups of mappings
$\Omega\to\Omega$ (subsemigroups of the \emph{full transformation semigroup}
$T(\Omega)$, where $\Omega:=\{1,\ldots ,n\}$). This area is especially
promising for two reasons. First, in a transformation
semigroup $S$, the units are the permutations; if there are any, they form a
\emph{permutation group} $G$ and we can take advantage of the very deep recent
results on them, chiefly the classification of finite simple groups (CFSG).
Secondly,  even if there are no units, we still have a group
to play with, the \emph{normaliser} of $S$ in $\Sym(\Omega)$, the set of
all permutations $g$ such that $g^{-1}Sg=S$.

The following result of Levi and McFadden \cite{lm} is the prototype for
results of this kind. Let $S_n$ and $T_n$ denote the symmetric group and full
transformation semigroup on $\Omega:=\{1,2,\ldots,n\}$.

\begin{thm}
Let $a\in T_n\setminus S_n$, and let $S$ be the semigroup generated by the
conjugates $g^{-1}ag$ for $g\in S_n$. Then
\begin{enumerate}
\item $S$ is idempotent-generated;
\item $S$ is regular;
\item $S=\langle a,S_n\rangle\setminus S_n$.
\end{enumerate}
\end{thm}

In other words, semigroups of this form, with normaliser $S_n$, have
\emph{very nice} properties!

Inspired by this result, we could formulate a general problem:

\begin{prob}\label{problem1}
\begin{enumerate}
\item
Given a semigroup property P, for which pairs $(a,G)$, with
$a\in T_n\setminus S_n$ and $G\le S_n$, does the semigroup
$\langle g^{-1}ag:g\in G\rangle$ have property P?
\item
Given a semigroup property P, for which pairs $(a,G)$ as above does the
semigroup $\langle a,G\rangle\setminus G$ have property P?
\item
For which pairs $(a,G)$ are the semigroups of the preceding parts equal?
\end{enumerate}
\end{prob}

The following portmanteau theorem lists some previously known results on this
problem. The first part is due to Levi \cite{levi96}, the other two to Ara\'ujo,
Mitchell and Schneider \cite{ArMiSc}.

\begin{thm}
\begin{enumerate}
\item
For any $a\in T_n\setminus S_n$ the semigroups
$\langle g^{-1}ag:g\in S_n\rangle$ and $\langle g^{-1}ag:g\in A_n\rangle$
are equal.
\item
$\langle g^{-1}ag:g\in G\rangle$ is
idempotent-generated for all $a\in T_n\setminus S_n$ if and only if $G=S_n$
or $G=A_n$ or $G$ is one of three specific groups of low degrees.
\item
$\langle g^{-1}ag:g\in G\rangle$ is regular
for all $a\in T_n\setminus S_n$ if and only if $G=S_n$ or $G=A_n$ or $G$ is
one of eight specific groups of low degrees.
\end{enumerate}
\end{thm}

Recently, we have obtained several extensions of these results. The first
theorem is proved in \cite{ArCa}.

\begin{thm}\label{arca}
Given $k$ with $1\le k\le n/2$, the following are equivalent for a subgroup
$G$ of $S_n$:
\begin{enumerate}
\item for all rank $k$ transformations $a$, $a$ is regular in
$\langle a,G\rangle$;
\item for all rank $k$ transformations $a$, $\langle a,G\rangle$ is
regular;
\item for all rank $k$ transformations $a$, $a$ is regular in
$\langle g^{-1}ag:g\in G\rangle$;
\item for all rank $k$ transformations $a$, $\langle g^{-1}ag:g\in G\rangle$
is regular.
\end{enumerate}
Moreover, we have a complete list of the possible groups $G$ with these
properties for $k\ge5$, and partial results for smaller values.
\end{thm}

It is worth pointing out that in the previous theorem the equivalence between
(a) and (c) is not new (it appears in \cite{lmm}).
Really surprising, and a great result that semigroups owe to the classification
of finite simple groups, are the equivalences between (a) and (b), and between
(c) and (d).

The four equivalent properties above translate into a transitivity property
of $G$ which we call the \emph{$k$-universal transversal property}, which we
will describe in the Subsection~\ref{kut}.

\spc

In the framework of Problem \ref{problem1}, let P be the following property:
the pair $(a,G)$, with $a\in T_n\setminus S_n$ and $G\le S_n$, satisfies
 $\langle a,G\rangle \setminus G = \langle a,S_n\rangle \setminus S_n$.

The classification of the pairs $(a,G)$ with this property poses a very interesting group theoretical problem. Recall that
the rank of a map $a\in T_n$ is $|\Omega a|$ and the kernel of $a$ is
$\ker(a):=\{(x,y)\in \Omega^2 \mid xa=ya\}$; by the usual correspondence
between equivalences and partitions, we can identify $\ker(a)$ with a
partition $\{A_1,\ldots , A_k\}$. Suppose $|\Omega|>2$ and we have a rank
$2$  map  $a\in T_n$. It is clear that $g^{-1}a\in \langle a,S_n\rangle$, for
all $g\in S_n$. In addition, if $\ker(a)=\{A_1,A_2\}$, then
$\ker(g^{-1}a)=\{A_1g,A_2g\}$. Therefore, in order to classify the groups
with property P above we need to find the groups $G$ such that
\begin{eqnarray}\label{lambda}
\{\{A_1,A_2\}g\mid g\in G\}&=&\{\{A_1,A_2\}g\mid g\in S_n\}.
\end{eqnarray}
If $|A_1|<|A_2|$, this is just $|A_1|$-homogeneity; but if these two sets have
the same size, the property is a little more subtle.

Extending this analysis to partitions with more than two parts, we see that the group-theoretic properties we need
to investigate are transitivity on ordered partitions of given shape (this
notion was introduced by Martin and Sagan~\cite{MS} under the name
\emph{partition-transitivity}) and the weaker notion of transitivity on
unordered partitions of given shape. This is done in Section~\ref{parthomog},
where we indicate the proof of the following theorem from~\cite{AnArCa}.

\begin{thm}
We have a complete list (in terms of the rank and kernel type of $a$) for
pairs $(a,G)$ for which $\langle a,G\rangle\setminus G=
\langle a,S_n\rangle\setminus S_n$.
\end{thm}

As we saw, the semigroups $\langle a,S_n\rangle\setminus S_n$
have very nice properties. In particular,
the questions of calculating their automorphisms and congruences, checking
for regularity, idempotent generation, etc., are all
settled. Therefore the same happens for the groups $G$ and maps
$a\in T_n\setminus S_n$ such that  $\langle a,G\rangle\setminus G=
\langle a,S_n\rangle\setminus S_n$, and all these pairs $(a,G)$ have been
classified.

Another long-standing open question was settled by the following theorem, from
\cite{ArCaMiNe}.

\begin{thm}\label{third}
The semigroups $\langle a,G\rangle\setminus G$ and
$\langle g^{-1}ag:g\in G\rangle$ are equal for all $a\in T_n\setminus S_n$
if and only if $G=S_n$, or $G=A_n$, or $G$ is the trivial group, or $G$ is
one of five specific groups.
\end{thm}

\begin{prob}
It would be good to have a more refined version of this where the hypothesis
refers only to all maps of rank $k$, or just a single map $a$.
\end{prob}

\subsection{Homogeneity and related properties}\label{homogeneity}

A permutation group $G$ on $\Omega$ is \emph{$k$-homogeneous} if it acts
transitively on the set of $k$-element subsets of $\Omega$, and is
\emph{$k$-transitive} if it acts transitively on the set of $k$-tuples of
distinct elements of $\Omega$.

It is clear that $k$-homogeneity is equivalent to $(n-k)$-homogeneity, where
$|\Omega|=n$; so we may assume that $k\le n/2$. It is also clear that
$k$-transitivity implies $k$-homogeneity.

We say that $G$ is \emph{set-transitive} if it is $k$-homogeneous for all
$k$ with $0\le k\le n$. The problem of determining the set-transitive groups
was posed by von Neumann and Morgenstern \cite{vNM} in the first edition of
their influential book on game theory. In the second edition, they refer to
an unpublished solution by Chevalley, but the first
published solution was by Beaumont and Peterson \cite{BP}. The set-transitive
groups are the symmetric and alternating groups, and four small exceptions
with degrees $5,6,9,9$.

In an elegant paper in 1965, Livingstone and Wagner~\cite{lw} showed:

\begin{thm}\label{lwthm}
Let $G$ be $k$-homogeneous, where $2\le k\le n/2$. Then
\begin{enumerate}
\item $G$ is $(k-1)$-homogeneous;
\item $G$ is $(k-1)$-transitive;
\item if $k\ge5$, then $G$ is $k$-transitive.
\end{enumerate}
\end{thm}

In particular, part (a) of this theorem is proved by a short argument using
character theory of the symmetric group. This can be translated into
combinatorics, and generalised to linear and affine groups: see
Kantor~\cite{kantor:inc}.

The $k$-homogeneous but not $k$-transitive groups for $k=2,3,4$ were
determined by Kantor~\cite{kantor:4homog,kantor:2homog}. All this was pre-CFSG.

The $k$-transitive groups for $k>1$ are known, but the classification uses
CFSG. Lists can be found in various references such as \cite{cam,dixon}.

\subsection{The $k$-universal transversal property}\label{kut}

Let $G\le S_n$, and $k$ an integer smaller than $n$.

The group $G$ has the \emph{$k$-universal transversal property}, or
\emph{$k$-ut} for short, if for every $k$-element subset $S$ of
$\{1,\ldots,n\}$ and every $k$-part partition $P$ of $\{1,\ldots,n\}$, there
exists $g\in G$ such that $Sg$ is a \emph{transversal} or \emph{section}
for $P$: that is, each part of $P$ intersects $Sg$ in a single point.

\begin{thm}
For $k\le n/2$, the following are equivalent for a permutation group $G\le S_n$:
\begin{enumerate}
\item for all $a\in T_n\setminus S_n$ with rank $k$, $a$ is regular in
$\langle a,G\rangle$;
\item $G$ has the $k$-universal transversal property.
\end{enumerate}
\end{thm}

In order to get the surprising equivalence (noted after Theorem \ref{arca})
 of ``$a$ is regular in $\langle a,G\rangle$''
and $``\langle a,G\rangle$ is regular'', we need to know that, for $k\le n/2$,
a group with the $k$-ut property also has the $(k-1)$-ut property. This fact,
the analogue of Theorem~\ref{lwthm}(a), is not at all obvious.

We go by way of a related property: $G$ is \emph{$(k-1,k)$-homogeneous} if,
given any two subsets $A$ and $B$ of $\{1,\ldots,n\}$ with $|A|=k-1$ and
$|B|=k$, there exists $g\in G$ with $Ag\subseteq B$.

Now the $k$-ut property implies $(k-1,k)$-homogeneity. (Take a partition with
$k$ parts, the singletons contained in $A$ and all the rest. If $Bg$ is a
transversal for this partition, then $Bg\supseteq A$, so $Ag^{-1}\subseteq B$.)

The bulk of the argument involves these groups. We show that, if
$3\le k\le(n-1)/2$ and $G$ is $(k-1,k)$-homogeneous, then either $G$ is
$(k-1)$-homogeneous, or $G$ is one of four small exceptions (with $k=3,4,5$ and
$n=2k-1$).

It is not too hard to show that such a group $G$ must be transitive, and
then primitive. Now careful consideration of the orbital graphs shows that
$G$ must be $2$-homogeneous, at which point we invoke the classification of
$2$-homogeneous groups (a consequence of CFSG).

One simple observation: if $G$ is $(k-1,k)$-homogeneous but not
$(k-1)$-homogeneous of degree $n$, then colour one $G$-orbit of $(k-1)$-sets
red and the others blue; by assumption, there is no monochromatic $k$-set,
so $n$ is bounded by the Ramsey number $R(k-1,k,2)$. The values
$R(2,3,2)=6$ and $R(3,4,2)=13$ are useful here; $R(4,5,2)$ is unknown, and
in any case too large for our purposes.

Now we return to considering the $k$-ut property.

First, we note that the
$2$-ut property says that every orbit on pairs contains a pair crossing
between parts of every $2$-partition; that is, every orbital graph is
connected. By Higman's Theorem, this is equivalent to primitivity.

For $2<k<n/2$, we know that the $k$-ut property lies between
$(k-1)$-homogeneity and $k$-homogeneity, with a few small exceptions. In
fact $k$-ut is equivalent to $k$-homogeneous for $k\ge6$; we classify all
the exceptions for $k=5$, but for $k=3$ and $k=4$ there are some groups
we are unable to resolve (affine, projective and Suzuki groups), which
pose interesting problems (see Problems~\ref{p:ut1} and \ref{p:ut2}).

For large $k$ we have:

\begin{thm}\label{1.10}
For $n/2<k<n$, the following are equivalent:
\begin{enumerate}
\item $G$ has the $k$-universal transversal property;
\item $G$ is $(k-1,k)$-homogeneous;
\item $G$ is $k$-homogeneous.
\end{enumerate}
\end{thm}

In the spirit of Livingstone and Wagner, we could ask:

\begin{prob}
Without using CFSG, show any or all of the following implications:
\begin{enumerate}
\item $k$-ut implies $(k-1)$-ut for $k\le n/2$;
\item $(k-1,k)$-homogeneous implies $(k-2,k-1)$-homogeneous for $k\le n/2$;

\item $k$-ut (or $(k-1,k)$-homogeneous) implies $(k-1)$-homogeneous for
$k\le n/2$.
\end{enumerate}
\end{prob}

\subsection{Partition transitivity and homogeneity}\label{parthomog}

Let $\lambda$ be a partition of $n$ (a non-increasing sequence of positive
integers with sum $n$). A partition of $\{1,\ldots,n\}$ is said to have
\emph{shape} $\lambda$ if the size of the $i$th part is the $i$th part of
$\lambda$.

The group $G$ is \emph{$\lambda$-transitive} if, given any two (ordered)
partitions of shape $\lambda$, there is an element of $G$ mapping each part
of the first to the corresponding part of the second. (This notion is due to
Martin and Sagan \cite{MS}.) Moreover, $G$ is
\emph{$\lambda$-homogeneous} if there is an element of $G$ mapping the
first partition to the second (but not necessarily respecting the order of
the parts).

Of course $\lambda$-transitivity implies $\lambda$-homogeneity, and the
converse is true if all parts of $\lambda$ are distinct.
If $\lambda=(n-t,1,\ldots,1)$, then $\lambda$-transitivity and
$\lambda$-homogeneity are equivalent to $t$-transitivity and $t$-homogeneity.

The connection with semigroups is given by the next result, from \cite{AnArCa}.
Let $G$ be a permutation group, and $a\in T_n\setminus S_n$, where $r$ is
the rank of $a$, and $\lambda$ the shape of the kernel partition.

\begin{thm}
For $G\le S_n$ and $a\in T_n\setminus S_n$, the following are equivalent:
\begin{enumerate}
\item $\langle a,G\rangle\setminus G=\langle a,S_n\rangle\setminus S_n$;
\item $G$ is $r$-homogeneous and $\lambda$-homogeneous.
\end{enumerate}
\end{thm}

So we need to know the $\lambda$-homogeneous groups. First, we consider
$\lambda$-transitive groups.

If the largest part of $\lambda$ is greater than $n/2$ (say $n-t$, where
$t<n/2$), then $G$ is $\lambda$-transitive if and only if it is
$t$-homogeneous
and the group $H$ induced on a $t$-set by its setwise stabiliser is
$\lambda'$-transitive, where $\lambda'$ is $\lambda$ with the part $n-t$
removed.

So if $G$ is $t$-transitive, then it is $\lambda$-transitive for all such
$\lambda$.

If $G$ is $t$-homogeneous but not $t$-transitive, then $t\le 4$, and
examination of the groups in Kantor's list gives the possible $\lambda'$
in each case.

So what remains is to show that, if $G$ is $\lambda$-transitive but not
$S_n$ or $A_n$, then $\lambda$ must have a part greater than $n/2$.

If $\lambda\ne(n),(n-1,1)$, then $G$ is primitive.

If $n\ge8$, then by \emph{Bertrand's Postulate}, there is a prime $p$ with
$n/2< p\le n-3$. If there is no part of $\lambda$ which is at least $p$, then
the number of partitions of shape $\lambda$ (and hence the order of $G$)
is divisible by $p$. A theorem of Jordan (see Wielandt~\cite{wie}, Theorem 13.9)
now shows that $G$ is symmetric or alternating.

The classification of $\lambda$-homogeneous but not $\lambda$-transitive
groups is a bit harder. We have to use
\begin{enumerate}
\item a little character theory to show that either $G$ fixes a point and is
transitive on the rest, or $G$ is transitive;
\item the argument using Bertrand's postulate and Jordan's theorem as before;
\item CFSG (to show that $G$ cannot be more than $5$-homogeneous if it is not
$S_n$ or $A_n$).
\end{enumerate}

The outcome is a complete list of such groups.

\subsection{Normalising groups}

We define a permutation group $G$ to be \emph{normalising} if
$\langle g^{-1}ag:g\in G\rangle=\langle a,G\rangle\setminus G$ for
all $a\in T_n\setminus S_n$.

The classification of normalising groups given by Theorem~%
\ref{third} is a little different; although permutation
group techniques are essential in the proof, we didn't find a simple
combinatorial condition on $G$ which is equivalent to this property.
We will not discuss it further here.

\section{Synchronization}

\subsection{Introduction}

In this section, we give a brief report on synchronization.

A (finite deterministic) \emph{automaton} consists of a finite set $\Omega$
of \emph{states} and a finite set of maps from $\Omega$ to $\Omega$ called
\emph{transitions}, which may be composed freely.

In other words, it is a transformation semigroup with a distinguished set of
generators.

An automaton is \emph{synchronizing} if there is a map of rank~$1$ (image of
size~$1$) in the semigroup. A word in the generators expressing such a map
is called a \emph{reset word}.

We will also call a transformation semigroup \emph{synchronizing} if it
contains an element of rank 1.

\begin{ex} This example has four (numbered) states, and two
transitions $A$ and $B$, shown as double and single lines respectively.

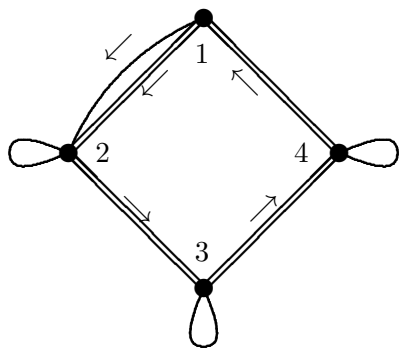
\begin{center}
\setlength{\unitlength}{1.2mm}
\begin{picture}(50,40)
\thicklines
\multiput(25,10)(0,30){2}{\circle*{2}}
\multiput(10,25)(30,0){2}{\circle*{2}}
\put(24,35){$1$}
\put(13,24){$2$}
\put(24,13){$3$}
\put(35,24){$4$}
\multiput(24.8,10.2)(-15,15){2}{\line(1,1){15}}
\multiput(25.2,9.8)(-15,15){2}{\line(1,1){15}}
\multiput(24.8,9.8)(15,15){2}{\line(-1,1){15}}
\multiput(25.2,10.2)(15,15){2}{\line(-1,1){15}}
\put(18,32){$\swarrow$}
\put(16,18){$\searrow$}
\put(30,18){$\nearrow$}
\put(28,32){$\nwarrow$}
\curve(10,25,16,34,25,40)
\put(14,36){$\swarrow$}
\curve(10,25,5,23.5,3.5,25,5,26.5,10,25)
\curve(25,10,23.5,5,25,3.5,26.5,5,25,10)
\curve(40,25,45,23.5,46.5,25,45,26.5,40,25)
\end{picture}
\end{center}
The reader can check easily that, irrespective of the starting state,
following the path $BAAABAAAB$ always ends in state $2$, and hence this is a
reset word of length~$9$. In fact, this is the shortest reset word.
\end{ex}

The \emph{\v{C}ern\'y Conjecture} asserts that if an $n$-state automaton is
synchronizing, then it has a reset word of length at most $(n-1)^2$. The above
example, with the square replaced by an $n$-gon, shows that this would be best
possible. The problem has been open for about 45 years. The best known bound
is cubic.

It is known that testing whether an automaton is synchronizing is in
\textsf{P}, but finding the length of the shortest reset word is
\textsf{NP}-hard.

\subsection{Graph homomorphisms and transformation semigroups}

All graphs here are undirected simple graphs (no loops or multiple edges).

A \emph{homomorphism} from a graph $X$ to a graph $Y$ is a map $f$ from the
vertex set of $X$ to the vertex set of $Y$ which carries edges to edges.
(We don't specify what happens to a non-edge; it may map to a non-edge, or
to an edge, or collapse to a vertex.) An \emph{endomorphism} of a graph $X$
is a homomorphism from $X$ to itself.

Let $K_r$ be the complete graph with $r$ vertices. The \emph{clique number}
$\omega(X)$ of $X$ is the size of the largest complete subgraph, and the
\emph{chromatic number} $\chi(X)$ is the least number of colours required for
a proper colouring of the vertices (adjacent vertices getting different
colours).

\begin{enumerate}
\item There is a homomorphism from $K_r$ to $X$ if and only if $\omega(X)\ge r$.
\item There is a homomorphism from $X$ to $K_r$ if and only if $\chi(X)\le r$.
\end{enumerate}

There are correspondences in both directions between graphs and transformation
semigroups (not quite functorial, or a Galois correspondence, sadly!)

First, any graph $X$ has an \emph{endomorphism semigroup} $\End(X)$.

In the other direction, given a transformation semigroup $S$ on $\Omega$, its
\emph{graph} $\Gr(S)$ has $\Omega$ as vertex set, two vertices $v$ and $w$
being joined if and only if there is no element of $S$ which maps $v$ and $w$
to the same place.

\begin{enumerate}
\item $\Gr(S)$ is complete if and only if $S\le S_n$;
\item $\Gr(S)$ is null if and only if $S$ is synchronizing;
\item $S\le\End(\Gr(S))$ for any $S\le T_n$;
\item $\omega(\Gr(S))=\chi(\Gr(S))$; this is equal to the minimum rank of an
element of $S$.
\end{enumerate}

Now the main theorem of this section describes the unique obstruction to
synchronization for a transformation semigroup.

\begin{thm}
A transformation semigroup $S$ on $\Omega$ is non-synchronizing if and only if
there is a non-null graph $X$ on the vertex set $\Omega$ with
$\omega(X)=\chi(X)$ and $S\le\End(X)$.
\end{thm}

In the reverse direction, the endomorphism semigroup of a non-null graph cannot
be synchronizing, since edges can't be collapsed. In the forward direction,
take $X=\Gr(S)$; there is some straightforward verification to do. (For details see \cite{ArCa2}.)

\subsection{Maps synchronized by groups}

Let $G\le S_n$ and $a\in T_n\setminus S_n$. We say that $G$ \emph{synchronizes}
$a$ if $\langle a,G\rangle$ is synchronizing.

By abuse of language, we say that $G$ is \emph{synchronizing} if it
synchronizes every element of $T_n\setminus S_n$.

Our main problem is to determine the synchronizing groups. From the theorem,
we see that $G$ is non-synchronizing if and only if there is a $G$-invariant
graph whose clique number and chromatic number are equal.

Rystsov \cite{rystsov} showed the following result, which implies that
synchronizing groups are necessarily primitive.

\begin{thm}
A permutation group $G$ of degree $n$ is primitive if and only if it
synchronizes every map of rank $n-1$.
\end{thm}

We give a brief sketch of the proof, to illustrate the graph endomorphism
technique. The backward implication is trivial; so suppose, for a
contradiction, that $G$ is primitive but fails to synchronize the map $a$ of
rank $n-1$. Then there are two points $x,y$ with $xa=ya$, and $a$ is bijective
on the remaining points. Choose a graph $X$ with $\langle G,a\rangle\le\End(X)$.
Note that $X$ is a regular graph. Since $a$ is an endomorphism, $x$ and $y$
are non-adjacent; so $a$ maps the neighbours of $x$ bijectively to the
neighbours of $xa$, and similarly the neighbours of $y$ to those of $ya$.
Since $xa=ya$, we see that $x$ and $y$ have the same neighbour set. Now
``same neighbour set'' is an equivalence relation preserved by $G$,
contradicting primitivity.

\spc

So a synchronizing group must be primitive.

We have recently improved this: a primitive group synchronizes every
map of rank $n-2$. The key tool in the proof is graph endomorphisms. Also, a
primitive group synchronizes every map of kernel type $(k,1,\ldots,1)$. For
both results, and further information, see \cite{ArCa2}.

Also, $G$ is synchronizing if and only if there is no $G$-invariant graph,
not complete or null, with clique number equal to chromatic number. For more
on this see \cite{ArnoldSteinberg,CK,neumann,JEP,rystsov,Tr,Tr07}. Thus,
a $2$-homogeneous group is synchronizing, and a synchronizing group is
primitive. For if $G$ is $2$-transitive, the only $G$-invariant graphs are
complete or null; and if $G$ is imprimitive, then it preserves a complete
multipartite graph.

Furthermore, a synchronizing group is \emph{basic} in the O'Nan--Scott
classification, that is, not contained in a wreath product with the product
action. (For non-basic primitive groups preserve Hamming graphs, which have
clique number equal to chromatic number.) By the O'Nan--Scott Theorem,
such a group is affine, diagonal or almost simple.

None of the above implications reverses. Indeed, there are non-synchronizing
basic groups of all three O'Nan--Scott types.

We are a long way from a classification of synchronizing groups. The attempts
to classify them lead to some interesting and difficult problems in
extremal combinatorics, finite geometry, computation, etc. But that is another
survey paper! We content ourselves here with a single result about an
important class of primitive groups, namely the classical symplectic,
orthogonal and unitary groups, acting on their associated polar spaces.
The implicit geometric problem has not been completely solved, despite decades
of work by finite geometers. We refer to Thas~\cite{jat} for a survey.

\begin{thm}
A classical group, acting on the points of its associated polar space, is
non-synchronizing if and only if the polar space possesses either an ovoid
and a spread, or a partition into ovoids.
\end{thm}

\subsection{A conjecture}

We regard the following as the biggest open problem in the area.
A map $a\in T_n$ is
\emph{non-uniform} if its kernel classes are not all of the same size.

\begin{conj}
A primitive permutation group synchronizes every non-uniform map.
\end{conj}

We have some partial results about this (see \cite{ABC,ArCa2}) but are far from a proof!

\section{Problems}\label{spro}

One of the goals of this paper is to provide a list of problems that
might help the interested reader involve himself in this fascinating topic. 
In addition to the problems included above, we collect here a number
of problems on the general interplay between properties of the group of
units and properties of the semigroup containing it.

We start by proposing a problem to experts in number theory. If this problem
can be solved, the results on $\agl(1,p)$, in \cite{ArCa}, will be dramatically
sharpened.
\begin{prob}\label{p:ut1}
Classify the prime numbers $p$ congruent to $11$ (mod $12$) such that for some
$c\in \GF(p)^*$ we have $|\langle -1,c,c-1\rangle|<p-1$.
\end{prob}

The primes less than $500$ with this property are $131$, $191$, $239$, $251$,
$311$, $419$, $431$, and $491$.

\begin{prob}\label{p:ut2}
Do the Suzuki groups $\Sz(q)$ have the $3$-universal transversal  property?

Classify the groups $G$ that have the $4$-ut property, when
$\psl(2,q) \le G \le\pgaml(2,q)$, with either $q$ prime (except
$\psl(2,q)$ for $q\equiv1$ (mod~$4$), which is not $3$-homogeneous), or $q=2^p$
for $p$ prime.
\end{prob}

A group $G\leq S_n$ has the $(n-1)$-universal transversal property if and only
if it is transitive. And $\langle a,G\rangle$ (for a rank $n-1$ map $a$) contains all the rank $n-1$ maps
of $T_n$ if and only if $G$ is $2$-homogeneous. In this last case
$\langle a,G\rangle$  is regular for all $a\in T_n$, because
$\langle a,G\rangle=\{b\in T_n\mid |\Omega b|\leq n-1\}\cup G$, and this
semigroup is well known to be regular.

\begin{prob}
Classify the groups $G\leq S_n$ such that $G$ together with any  rank $n-k$
map, where $k\leq 5$,  generate a regular semigroup. We already know that
such $G$ must be $k$-homogeneous; so we know which groups to look at (see Theorem \ref{1.10}).
\end{prob}

The difficulty here (when rank $k>\lfloor \frac{n+1}{2}\rfloor$) is that a
$k$-homogenous group is not necessarily $(k-1)$-homogenous. Therefore a rank
$k$ map $a\in T_n$ might be regular in $\langle a, G\rangle$, but we are not
sure that there exists  $g\in G$ such that $\rank(bgb)=\rank(b)$,
for $b\in  \langle a, G\rangle$ such that $\rank(b)<\rank(a)$.

It is clear that if $\langle a ,G\rangle\setminus G$ is idempotent generated, for all rank $k$ transformation $a\in T_n\setminus S_n$, then $G$ has the $k$-ut property (see \cite{ArMiSc}).
 
\begin{prob}
Classify the groups $G\leq S_n$ such that $\langle a ,G\rangle\setminus G$ is
idempotent generated, for all  rank $k$ maps, where $k\leq n/2$. Even if the
classification of the groups with the $k$-ut property is only almost finished
(Problem \ref{p:ut2} is the missing part), it might be possible to settle
the idempotent generation problem.
\end{prob}

\begin{prob}\label{top}
The most general problem that has to be handled is the classification of pairs
$(a,G)$, where $a\in T_n$ and $G\leq S_n$, such that $\langle a,G\rangle$ is a
regular semigroup.
\end{prob}

When investigating $(k-1)$-homogenous groups without the $k$-universal
transversal property ($k$-ut property), it was common that some of the orbits
on the $k$-sets have transversals for all the partitions. Therefore the
following definition is natural.

A group $G\leq S_n$ is said to have the weak $k$-ut property if there exists
a $k$-set $S\subseteq \Omega$ such that the orbit of $S$ under $G$ contains a
section for all $k$-partitions. Such a set is called a $G$-universal
transversal set. A solution to the following problem would have important
consequences in semigroup theory.

\begin{prob}\label{five}
Classify the groups with the weak $k$-ut property; in addition, for each
of them, classify their $G$-universal transversal sets.
\end{prob}

In McAlister's celebrated paper \cite{mcalister} it is proved that, if
$e^2=e\in T_n$ is a rank $n-1$ idempotent, then $\langle G,e\rangle$ is
regular for all groups $G\leq S_n$. In addition, assuming  that
$\{\alpha,\beta\}$ is the non-singleton kernel class of $e$ and
$\alpha e=\beta$,  if $\alpha$ and $\beta$ are not in the same orbit under
$G$, then   $\langle e,G\rangle$ is an orthodox semigroup (that is, the
idempotents form a subsemigroup); and $\langle e,G\rangle$ is inverse if and
only if $\alpha$ and $\beta$ are not in the same orbit under $G$ and the
stabilizer of $\alpha$ is contained in the stabilizer of $\beta$.

\begin{prob}
Classify the groups $G\leq S_n$ that together with any idempotent [rank $k$
idempotent] generate a regular [orthodox, inverse] semigroup.

Classify the pairs $(G,a)$, with $a\in T_n$ and $G\leq S_n$, such that
$\langle e,G\rangle$ is inverse [orthodox].
\end{prob}

The theorems and problems in this paper admit linear versions that are
interesting for experts in groups and semigroups, but also to experts in
linear algebra and matrix theory.

\begin{prob}
Prove (or disprove) that if $G\leq \gl(n,q)$  such that for all singular
matrix $a$ there exists $g\in G$ with  $\rank(a)=\rank(aga)$, then $G$
contains the special linear group.
\end{prob}

For $n=2$ and for $n=3$, this condition is equivalent to irreducibility of $G$.
But we conjecture that, for sufficiently large $n$, it implies that $G$
contains the special linear group.

\begin{prob}
Classify the groups $G\leq \gl(n,q)$  such that for all rank $k$ (for a given
$k$) singular matrix $a$ we have that $a$ is regular in $\langle G,a\rangle$
[the semigroup $\langle G,a\rangle$ is regular].
\end{prob}

To handle this problem it is useful to keep in mind the following results.
Kantor~\cite{kantor:inc} proved that if a subgroup of $\pgaml(d,q)$ acts
transitively on $k$-dimensional subspaces, then it acts transitively on
$l$-dimensional subspaces for all $l\le k$ such that $k+l\le n$;
in~\cite{kantor:line}, he showed that subgroups transitive on $2$-dimensional
subspaces are $2$-transitive on the $1$-dimensional subspaces with the single
exception of a subgroup of $\pgl(5,2)$ of order $31\cdot5$; and, with the
second author~\cite{cameron-kantor}, he showed that such groups must contain
$\psl(d,q)$ with the single exception of the alternating group $A_7$ inside
$\pgl(4,2)\cong A_8$. Also Hering \cite{He74,He85} and Liebeck \cite{Li86},
using CFSG,
classified the subgroups of $\pgl(d,p)$ which are transitive on $1$-spaces.

\spc

Regarding synchronization, the most important question (in our opinion) is
the following conjecture, stated earlier.

\begin{prob}
Is it true that every primitive group of permutations of a finite set $\Omega$
synchronizes every non-uniform transformation on $\Omega$?
\end{prob}

Assuming the previous question has an affirmative answer (as we believe), an
intermediate step in order to prove it would be to solve the following set of
connected problems:

\begin{prob}
\begin{enumerate}
\item Prove that every map of rank $n-3$, with non-uniform kernel, is
synchronized by a primitive group. This is known for idempotent maps (see
\cite{ArCa2}).
 \item Prove that a primitive group synchronizes every non-uniform map of rank
$5$.
 \item Prove that if in $S=\langle f,G\rangle$ there  is a map of minimal rank
$r>1$,  there can be no map in $S$ with rank $r+2$.
\end{enumerate}
\end{prob}

The next class of groups lies strictly between primitive and synchronizing.

\begin{prob}
Is it possible to classify the primitive groups which synchronize
every rank $3$ map?
\end{prob}
		
Note that there are primitive groups that do not synchronize a
rank $3$ map (see \cite{neumann}). And there are non-synchronizing groups
which synchronize every rank $3$ map. Take for example $\pgl(2,7)$ of degree
$28$; this group is non-synchronizing, but synchronizes every rank $3$
map, since $28$ is not divisible by $3$.

\spc

There are very fast polynomial-time algorithms to decide if a given set of
permutations generates a primitive group, or a $2$-transitive group.

\begin{prob}
Find an efficient algorithm to decide if a given set of permutations generates
a synchronizing group.
\end{prob}

It would be quite remarkable if such an algorithm exists; as we saw, it would
in particular resolve questions about ovoids and spreads in certain polar
spaces (among other things).

\spc

There are a number of natural problems related to $\lambda$-homogeneity.

\begin{prob}
Let $H\leq \sym$ be a $2$-transitive group. Classify the pairs $(a,G)$, where
$a\in \sym $ and $G\leq \sym$, such that  $\langle a,G\rangle  = H$.
\end{prob}

\begin{prob}		
Let $G\leq \sym$ be a $2$-transitive group. (The list of those groups is
available in \cite{cam,dixon}.) For every $a\in \trans$ describe the structure
of $\langle G,a\rangle\setminus G$.
In particular (where $G$ is a $2$-transitive group and $a\in \trans$):
\begin{enumerate}
\item  classify all the pairs $(a,G)$ such that $\langle a,G\rangle $ is a
regular semigroup (that is, for all $x\in  \langle a,G\rangle$ there exists
$y \in \langle a,G\rangle$ such that $x=xyx$);
\item classify all the pairs $(a,G)$ such that $\langle a,G\rangle\setminus G$
is generated by its idempotents;
\item classify all the pairs $(a,G)$ such that
$\langle a,G\rangle\setminus G =\langle g^{-1}ag\mid g\in G\rangle $;
\item describe the automorphisms, congruences, principal right, left and
two-sided ideals of the semigroups $\langle a, G\rangle$ (when $G$ is a
$2$-transitive group).
\end{enumerate}
\end{prob}

\begin{prob}
For each $2$-transitive group $G$ classify the $G$-pairs, that is, the pairs
$(a,H)$ such that  $H\leq \sym$, $a\in\trans$  and
$\langle a,G\rangle\setminus G=\langle a,H\rangle\setminus H$.
\end{prob}

\begin{prob}
Let $V$ be a finite dimension vector space. A pair $(a,G)$, where $a$ is a
singular endomorphism of $V$ and  $G\leq\Aut(V)$, is said to be an
$\Aut(V)$-pair if
$$\langle a,G\rangle\setminus G=\langle a,\Aut(V)\rangle\setminus\Aut(V).$$

Classify  the $\Aut(V)$-pairs.
\end{prob}

\begin{prob}\label{11}
Formulate and prove analogues of the results in this paper, but for semigroups
of linear maps on a vector space.
\end{prob}

\begin{prob}
Solve the analogue of Problem \ref{11} for independence algebras (for
definitions and fundamental results see \cite{ArEdGi,arfo,cameronSz,F1,gould}).
\end{prob}

\textbf{Acknowledgment}
The first  author was partially supported by
Pest-OE/MAT/UI0143/2011 of Centro de Algebra da Universidade de Lisboa, and
by FCT and PIDDAC through the project PTDC/MAT/101993/2008.

\end{document}